\documentclass[11pt,reqno]{amsart}
\usepackage{amsfonts,amssymb,amsmath}
\usepackage{color}

\DeclareMathOperator{\bmax}{\mathbf{max}}

\DeclareMathOperator{\acyclic}{\mathrm{acyclic}}

\newtheorem{theorem}{Theorem}

\newtheorem{proposition}[theorem]{Proposition}

\numberwithin{equation}{section} \numberwithin{theorem}{section}

\begin{document}

\title{Pattern Recognition on Oriented Matroids:
Decompositions of Topes,\\ and Dehn--Sommerville Type Relations}

\author{Andrey O. Matveev}
\email{andrey.o.matveev@gmail.com}



\begin{abstract}
If $\mathrm{V}(\boldsymbol{R})$ is the vertex set of a symmetric cycle $\boldsymbol{R}$ in the tope graph of a simple oriented matroid $\mathcal{M}$, then for any tope $T$ of $\mathcal{M}$ there exists a unique inclusion-minimal subset $\boldsymbol{Q}(T,\boldsymbol{R})$ of $\mathrm{V}(\boldsymbol{R})$ such that $T$ is the sum of the topes of $\boldsymbol{Q}(T,\boldsymbol{R})$.

If $|\boldsymbol{Q}(T,\boldsymbol{R})|\geq 5$, then the decomposition $\boldsymbol{Q}(T,\boldsymbol{R})$ of the tope $T$ with respect to the symmetric cycle $\boldsymbol{R}$ satisfies certain Dehn--Sommerville type relations.
\end{abstract}

\maketitle

\pagestyle{myheadings}

\markboth{A.O.~MATVEEV}{PATTERN RECOGNITION ON ORIENTED MATROIDS}

\thispagestyle{empty}


\section{Introduction}

Let $\mathcal{A}:=(E_t,\mathcal{L})=(E_t,\mathcal{T})$ be a {\em simple\/} (i.e., with no {\em loops}, {\em parallel elements\/} or {\em antiparallel elements})
{\em acyclic\/} oriented matroid on the ground set~$E_t:=\{1,2,\ldots,t\}$, with set of covectors~$\mathcal{L}$ and with set of topes~$\mathcal{T}\subseteq\{1,-1\}^{E_t}$, of rank~$r\geq 3$. Let $\boldsymbol{R}$ be a {\em symmetric $2t$-cycle}, with its vertex set~$\mathrm{V}(\boldsymbol{R})=-\mathrm{V}(\boldsymbol{R})$, in the tope graph $\mathcal{T}(\mathcal{L}(\mathcal{A}))$ of~$\mathcal{A}$
such that~$(1,1,\ldots,1)=:\mathrm{T}^{(+)}\not\in\mathrm{V}(\boldsymbol{R})$.

Recall that the vertex set $\mathrm{V}(\boldsymbol{R})$ of the cycle $\boldsymbol{R}$ is the set of topes of a~{\em rank~$2$ oriented matroid\/} denoted $\mathcal{A}_{\boldsymbol{R}}$, see e.g.~\cite[Rem.~1.7]{M-DG-I}, and thus by~\cite[Example~7.1.7]{BLSWZ} we can consider the representation of $\mathcal{A}_{\boldsymbol{R}}$ by a certain central line arrangement
\begin{equation}
\label{prchap:9:eq:3}
\{\mathbf{x}\in\mathbb{R}^2\colon \langle\boldsymbol{a}_e,\mathbf{x}\rangle=0,\; e\in E_t\}
\end{equation}
in the plane~$\mathbb{R}^2$, with the set of corresponding normal vectors
\begin{equation}
\label{prchap:9:eq:7}
\boldsymbol{A}:=\{\boldsymbol{a}_e\colon e\in E_t\}\; ,
\end{equation}
and with the standard scalar product $\langle\cdot,\cdot\rangle$.

For any tope $T$ of the oriented matroid~$\mathcal{A}$ there exists a {\em unique inclusion-minimal subset\/} $\boldsymbol{Q}(T,\boldsymbol{R})\subset\mathrm{V}(\boldsymbol{R})$ such that
\begin{equation*}
\sum_{Q\in\boldsymbol{Q}(T,\boldsymbol{R})}Q=T\; ,
\end{equation*}
see~\cite[Sect.~11.1]{M-DG-I}.

If the tope graph $\mathcal{T}(\mathcal{L}(\mathcal{A}))$ is the {\em hypercube graph\/} $\boldsymbol{H}(t,2)$ on~$2^t$ vertices, then for any {\em odd\/} integer $j$, $1\leq j\leq t$, by~\cite[Th.~13.6]{M-DG-I} there are precisely~$2\tbinom{t}{j}$ vertices $T$ of~$\boldsymbol{H}(t,2)$ such that~$|\boldsymbol{Q}(T,\boldsymbol{R})|=j$.

Let $\bmax^+\mathrm{V}(\boldsymbol{R})$ denote the subset of topes of $\mathrm{V}(\boldsymbol{R})$ with {\em inclusion-maximal positive parts}. For the {\em positive tope\/} $\mathrm{T}^{(+)}$ of $\mathcal{A}$, by~\cite[Prop.~1.12(ii)]{M-DG-I} we have
\begin{equation*}
\boldsymbol{Q}(\mathrm{T}^{(+)},\boldsymbol{R})=\bmax^+\mathrm{V}(\boldsymbol{R})\; .
\end{equation*}

If
\begin{equation*}
|\boldsymbol{Q}(\mathrm{T}^{(+)},\boldsymbol{R})|\geq 5\; ,
\end{equation*}
then for the set of vectors~(\ref{prchap:9:eq:7}) we have
\begin{equation}
\label{prchap:9:eq:1}
|\{\boldsymbol{a}\in\boldsymbol{A}\colon \boldsymbol{a}\in\mathbf{C}_{>}\}|\geq 2\; ,
\end{equation}
for any {\em open half-space\/} $\mathbf{C}_{>}\subset \mathbb{R}^2$ bounded by a $1$-dimensional subspace of~$\mathbb{R}^2$; see e.g.~\cite[Proof of Prop.~2.33]{G-DG}.

Denote by $\nu_j$ the number of {\em feasible\/} subsystems, of cardinality $j$, of the {\em infeasible system\/} of {\em homogeneous strict linear inequalities}
\begin{equation}
\label{prchap:9:eq:2}
\bigl\{\langle\boldsymbol{a}_e,\mathbf{x}\rangle>0\colon\mathbf{x}\in\mathbb{R}^2,\; \boldsymbol{a}_e\in\boldsymbol{A}\bigr\}
\end{equation}
associated with the arrangement~(\ref{prchap:9:eq:3}). Since
the condition~(\ref{prchap:9:eq:1}) holds, the quantities~$\nu_j$ satisfy the relations (where~$\mathrm{x}$ is a formal variable)
\begin{equation}
\label{prchap:9:eq:4}
\begin{cases}
\nu_j=\tbinom{t}{j},\text{\rm\ if $0\leq j\leq 2$}\; ,\\ \nu_{t-1}=\nu_t=0\; ,\\
\sum_{j=3}^t\left(\tbinom{t}{j}\!-\!\nu_j\right)(\mathrm{x}\!-\!1)^{t-j}\!=\!
-\sum_{j=3}^t(-1)^j\left(\tbinom{t}{j}\!-\!\nu_j\right)\mathrm{x}^{t-j}
\end{cases}
\end{equation}
called the {\em Dehn--Sommerville equations for the feasible subsystems of the system\/} (\ref{prchap:9:eq:2}); see e.g.~\cite[Prop.~3.53]{G-DG}.

In particular,
\begin{itemize}
\item[$\bullet$]
if $t=5$, then
\begin{equation*}
\nu_{3}=\tbinom{t}{2}-t\; ;
\end{equation*}
\item[$\bullet$]
if $t=6$, then
\begin{align*}
\nu_{3}&=\tbinom{t}{3}-2t+4\; ,\\ \nu_{4}&=\tbinom{t}{2}-3t+6\; ;
\end{align*}
\item[$\bullet$]
if $t=7$, then
\begin{align*}
\nu_{4}&=2\nu_{3}-2\tbinom{t}{4}+\tbinom{t}{3}\; ,\\
\nu_{5}&=\nu_{3}-\tbinom{t}{4}+\tbinom{t}{2}-t\; .
\end{align*}
\end{itemize}

For any $j$, $3\leq j\leq t-2$, by~\cite[Cor.~3.55(i)]{G-DG} we have
\begin{equation*}
\tbinom{t}{j}\!-\!\nu_j=-\sum_{i=3}^j(-1)^i\tbinom{t-i}{j-i}\Bigl(\tbinom{t}{i}\!-\!\nu_i\Bigr)\; .
\end{equation*}

Recall also that
\begin{equation*}
\sum_{j=1}^{t-2}(-1)^j\nu_j=0\; .
\end{equation*}

\section{Decompositions of topes with respect to symmetric cycles in the tope graph, and Dehn--Sommerville type relations}

Let $\mathcal{M}:=(E_t,\mathcal{L})=(E_t,\mathcal{T})$ be an arbitrary {\em simple\/}
oriented matroid of rank~$r\geq 3$, and let $\boldsymbol{R}$ be a {\em symmetric $2t$-cycle\/} in the tope graph $\mathcal{T}(\mathcal{L}(\mathcal{M}))$ of~$\mathcal{M}$.

Given a tope $T\in\mathcal{T}$ of $\mathcal{M}$, consider the corresponding {\em decomposition
\begin{equation*}
T=\sum_{Q\in\boldsymbol{Q}(T,\boldsymbol{R})}Q
\end{equation*}
of\/ $T$ with respect to the cycle~$\boldsymbol{R}$}, for a {\em unique inclusion-minimal subset\/}~$\boldsymbol{Q}(T,\boldsymbol{R})$ of the vertex set $\mathrm{V}(\boldsymbol{R})$ of the cycle~$\boldsymbol{R}$. Suppose that
\begin{equation*}
|\boldsymbol{Q}(T,\boldsymbol{R})|\geq 5\; ,
\end{equation*}
and consider the {\em totally cyclic\/} oriented matroid
\begin{equation*}
\mathcal{Y}_{T,\boldsymbol{R}}:={}_{-T^-}(\mathcal{M}_{\boldsymbol{R}})
\end{equation*}
obtained from the rank~$2$ oriented matroid $\mathcal{M}_{\boldsymbol{R}}:=(E_t,\mathrm{V}(\boldsymbol{R}))$ with the set of topes~$\mathrm{V}(\boldsymbol{R})$
by {\em reorientation\/} of~$\mathcal{M}_{\boldsymbol{R}}$ on the {\em negative part\/} $T^{-}$ of the tope~$T$.

Associate with the abstract simplicial {\em complex}
\begin{equation}
\label{prchap:9:eq:5}
\Delta:=\Delta_{\acyclic}(\mathcal{Y}_{T,\boldsymbol{R}})
\end{equation}
of {\em acyclic subsets\/} of the ground set of the oriented matroid $\mathcal{Y}_{T,\boldsymbol{R}}$ its~``{\em long}'' {\em $f$-vector\/}
\begin{equation*}
\boldsymbol{f}(\Delta;t):=\bigl(f_0(\Delta;t),f_1(\Delta;t),\ldots,f_t(\Delta;t)\bigr)\in\mathbb{N}^{t+1}
\end{equation*}
defined by
\begin{equation*}
f_j(\Delta;t):=\#\{F\in\Delta\colon |F|=j\}\; ,\ \ \ 0\leq j\leq t\; .
\end{equation*}
In view of~(\ref{prchap:9:eq:4}), we have
\begin{equation*}
\begin{cases}
f_j(\Delta;t)=\tbinom{t}{j},\text{\rm\ if $0\leq j\leq 2$}\; ,\\ f_{t-1}(\Delta;t)=f_t(\Delta;t)=0\; ,\\
\sum_{j=3}^t\left(\tbinom{t}{j}\!-\!f_j(\Delta;t)\right)(\mathrm{x}\!-\!1)^{t-j}\!=\!
-\sum_{j=3}^t(-1)^j\left(\tbinom{t}{j}\!-\!f_j(\Delta;t)\right)\mathrm{x}^{t-j}\; .
\end{cases}
\end{equation*}

Let us return to the symmetric cycle~$\boldsymbol{R}$, and associate with the tope $T$ the abstract simplicial complex $\boldsymbol{\Lambda}:=\boldsymbol{\Lambda}(T,\boldsymbol{R})$ with the facet family
\begin{equation}
\label{prchap:9:eq:6}
\{E_t-\mathbf{S}(T,Q)\colon Q\in\boldsymbol{Q}(T,\boldsymbol{R})\}\; ,
\end{equation}
where $\mathbf{S}(T,Q):=\{e\in E_t\colon T(e)\neq Q(e)\}$ is the {\em separation set} of the topes~$T$ and~$Q$. By construction, the complex $\boldsymbol{\Lambda}$ and the complex~$\Delta$ defined by~(\ref{prchap:9:eq:5}) coincide, and we come to the following result:

\begin{proposition}
Let $\boldsymbol{R}$ be a symmetric cycle in the tope graph of a simple oriented matroid~$\mathcal{M}:=(E_t,\mathcal{T})$.

Let $T\in\mathcal{T}$ be a tope of $\mathcal{M}$ such that for the unique inclusion-minimal subset of topes $\boldsymbol{Q}(T,\boldsymbol{R})\subset\mathrm{V}(\boldsymbol{R})$ with the property
\begin{equation*}
\sum_{Q\in\boldsymbol{Q}(T,\boldsymbol{R})}Q=T
\end{equation*}
we have
\begin{equation*}
|\boldsymbol{Q}(T,\boldsymbol{R})|\geq 5\; .
\end{equation*}

\noindent{\rm(i)} The components of the long $f$-vector $\boldsymbol{f}(\boldsymbol{\Lambda};t)$ of the complex
\begin{equation*}
\begin{split}
\boldsymbol{\Lambda}:\!&=\boldsymbol{\Lambda}(T,\boldsymbol{R})\\
&=\boldsymbol{\Lambda}(-T,\boldsymbol{R})
\end{split}
\end{equation*}
whose family of facets is defined by~{\rm(\ref{prchap:9:eq:6})} satisfy the {\em Dehn--Sommerville} type {\em relations}
\begin{equation*}
\begin{cases}
f_j(\boldsymbol{\Lambda};t)=\tbinom{t}{j},\text{\rm\ if $0\leq j\leq 2$}\; ,\\ f_{t-1}(\boldsymbol{\Lambda};t)=f_t(\boldsymbol{\Lambda};t)=0\; ,\\
\sum_{j=3}^t\left(\tbinom{t}{j}\!-\!f_j(\boldsymbol{\Lambda};t)\right)(\mathrm{x}\!-\!1)^{t-j}\!=\!
-\sum_{j=3}^t(-1)^j\left(\tbinom{t}{j}\!-\!f_j(\boldsymbol{\Lambda};t)\right)\mathrm{x}^{t-j}\; .
\end{cases}
\end{equation*}

In particular,
\begin{itemize}
\item[$\bullet$]
if $t=5$, then
\begin{equation*}
f_3(\boldsymbol{\Lambda};t)=\tbinom{t}{2}-t=5\; ;
\end{equation*}
\item[$\bullet$]
if $t=6$, then
\begin{align*}
f_3(\boldsymbol{\Lambda};t)&=\tbinom{t}{3}-2t+4=12\; ,\\ f_4(\boldsymbol{\Lambda};t)&=\tbinom{t}{2}-3t+6=3\; ;
\end{align*}
\item[$\bullet$]
if $t=7$, then
\begin{equation*}
\renewcommand{\arraystretch}{1.4}
\begin{array}{cclcl}
f_3(\boldsymbol{\Lambda};t)&=&\tfrac{1}{2}f_4(\boldsymbol{\Lambda};t)+\tbinom{t}{4}-\tfrac{1}{2}\tbinom{t}{3}&=&\tfrac{1}{2}\bigl(f_4(\boldsymbol{\Lambda};t)+35\bigr)\\
&=&f_5(\boldsymbol{\Lambda};t)+\tbinom{t}{4}-\tbinom{t}{2}+t&=&f_5(\boldsymbol{\Lambda};t)+21\; ,\\
f_4(\boldsymbol{\Lambda};t)&=&2f_3(\boldsymbol{\Lambda};t)-2\tbinom{t}{4}+\tbinom{t}{3}&=&2f_3(\boldsymbol{\Lambda};t)-35\\
&=&2f_5(\boldsymbol{\Lambda};t)+\tbinom{t}{3}-2\tbinom{t}{2}+2t&=&2f_5(\boldsymbol{\Lambda};t)+7\; ,\\
f_5(\boldsymbol{\Lambda};t)&=&f_3(\boldsymbol{\Lambda};t)-\tbinom{t}{4}+\tbinom{t}{2}-t&=&f_3(\boldsymbol{\Lambda};t)-21\\
&=&\tfrac{1}{2}f_4(\boldsymbol{\Lambda};t)-\tfrac{1}{2}\tbinom{t}{3}+\tbinom{t}{2}-t&=&\tfrac{1}{2}\bigl(f_4(\boldsymbol{\Lambda};t)-35\bigr)+14\; ;
\end{array}
\renewcommand{\arraystretch}{1.0}
\end{equation*}
the quantity $f_4(\boldsymbol{\Lambda};7)$ is {\em odd}.
\end{itemize}

\noindent{\rm(ii)} For any $j$, $3\leq j\leq t-2$, we have
\begin{equation*}
\tbinom{t}{j}\!-\!f_j(\boldsymbol{\Lambda};t)=-\sum_{i=3}^j(-1)^i\tbinom{t-i}{j-i}\Bigl(\tbinom{t}{i}\!-\!f_i(\boldsymbol{\Lambda};t)\Bigr)\; .
\end{equation*}

\noindent{\rm(iii)} We have
\begin{equation*}
\sum_{j=1}^{t-2}(-1)^j f_j(\boldsymbol{\Lambda};t)=0\; .
\end{equation*}
\end{proposition}

\vspace{3mm}


\begin{thebibliography}{10}
\bibitem{BLSWZ}
{\em Bj\"{o}rner A.}, {\em Las~Vergnas M.}, {\em Sturmfels B.}, {\em White N.}, {\em Ziegler G.M.} Oriented matroids. Second
edition. Encyclopedia of Mathematics, 46. -- Cambridge: Cambridge University Press, 1999.

\bibitem{G-DG}
{\em Gainanov D.N.} Graphs for pattern recognition. Infeasible systems of linear inequalities. -- Berlin: De Gruyter, 2016.

\bibitem{M-DG-I}
{\em Matveev A.O.} Pattern recognition on oriented matroids. -- Berlin: De Gruyter, 2017.
\end{thebibliography}
\end{document}